\documentclass[a4paper,12pt]{amsart}
\usepackage{amsfonts}
\usepackage{amssymb}
\usepackage{ifthen}
\usepackage{graphicx}
\nonstopmode \numberwithin{equation}{section}
\usepackage[usenames]{color}
\newtheorem{thm}{Theorem}[section]
\newtheorem{lem}{Lemma}[section]
\newtheorem{cor}[thm]{Corollary}
\newtheorem{prop}[thm]{Proposition}

\newtheorem{step}{Step}[section]

\theoremstyle{definition}
\newtheorem{mlem}{Main lemma}[section]
\newtheorem{assertion}{Assertion}[section]
\newtheorem{cl}{Claim}[section]
\newtheorem{ca}{Case}[section]
\newtheorem{sca}{Subcase}[section]
\newtheorem{scl}{Subclaim}[section]
\newtheorem{conj}[thm]{Conjecture}
\newtheorem{fact}{Fact}[section]
\newtheorem{defn}[thm]{Definition}
\newtheorem{op}[thm]{Open Problem}
\newtheorem{ques}[thm]{Question}
\newtheorem{rem}[thm]{Remark}
\newtheorem{exam}[thm]{Example}

\numberwithin{equation}{section}

\newcounter {own}
\def\theown {\thesection       .\arabic{own}}

\newenvironment{pf}[1][]{%
 \vskip 3mm
 \noindent
 \ifthenelse{\equal{#1}{}}%
  {{\slshape Proof. }}%
  {{\slshape #1.} }%
 }%
{\qed\bigskip}

\newcounter{alphabet}
\newcounter{tmp}
\newenvironment{Thm}[1][]{\refstepcounter{alphabet}%
\bigskip%
\noindent%
{\bf Theorem \Alph{alphabet}}%
\ifthenelse{\equal{#1}{}}{}{ (#1)}%
{\bf .} \itshape}{\vskip 8pt}

\makeatletter
\newcommand{\Ref}[1]{\@ifundefined{r@#1}{}{\setcounter{tmp}{\ref{#1}}\Alph{tmp}}}
\makeatother

\newenvironment{Lem}[1][]{\refstepcounter{alphabet}%
\bigskip%
\noindent%
{\bf Lemma \Alph{alphabet}}%
{\bf .} \itshape}{\vskip 8pt}

\newcounter{alphabet2}

\def\be{\begin{equation}}
\def\ee{\end{equation}}

\newcommand{\ben}{\begin{enumerate}}
\newcommand{\een}{\end{enumerate}}

\newcommand{\blem}{\begin{lem}}
\newcommand{\elem}{\end{lem}}
\newcommand{\bthm}{\begin{thm}}
\newcommand{\ethm}{\end{thm}}
\newcommand{\bcor}{\begin{cor}}
\newcommand{\ecor}{\end{cor}}
\newcommand{\beg}{\begin{exam}}
\newcommand{\eeg}{\end{exam}}
\newcommand{\begs}{\begin{examples}}
\newcommand{\eegs}{\end{examples}}
\newcommand{\bdefe}{\begin{defn}}
\newcommand{\edefe}{\end{defn}}
\newcommand{\bprob}{\begin{prob}}
\newcommand{\eprob}{\end{prob}}
\newcommand{\bques}{\begin{ques}}
\newcommand{\eques}{\end{ques}}
\newcommand{\bei}{\begin{itemize}}
\newcommand{\eei}{\end{itemize}}
\newcommand{\bcon}{\begin{conj}}
\newcommand{\econ}{\end{conj}}
\newcommand{\bop}{\begin{op}}
\newcommand{\eop}{\end{op}}

\newcommand{\bas}{\begin{assertion}}
\newcommand{\eas}{\end{assertion}}

\newcommand{\bfa}{\begin{fact}}
\newcommand{\efa}{\end{fact}}

\newcommand{\bca}{\begin{ca}}
\newcommand{\eca}{\end{ca}}

\newcommand{\bst}{\begin{step}}
\newcommand{\est}{\end{step}}

\newcommand{\bsca}{\begin{sca}}
\newcommand{\esca}{\end{sca}}

\newcommand{\bcl}{\begin{cl}}
\newcommand{\ecl}{\end{cl}}

\newcommand{\bmlem}{\begin{mlem}}
\newcommand{\emlem}{\end{mlem}}

\newcommand{\bscl}{\begin{scl}}
\newcommand{\escl}{\end{scl}}

\newcommand{\bcons}{\begin{conjs}}
\newcommand{\econs}{\end{conjs}}

\newcommand{\bprop}{\begin{prop}}
\newcommand{\eprop}{\end{prop}}

\newcommand{\br}{\begin{rem}}
\newcommand{\er}{\end{rem}}
\newcommand{\brs}{\begin{rems}}
\newcommand{\ers}{\end{rems}}
\newcommand{\bo}{\begin{obser}}
\newcommand{\eo}{\end{obser}}
\newcommand{\bos}{\begin{obsers}}
\newcommand{\eos}{\end{obsers}}
\newcommand{\bpf}{\begin{pf}}
\newcommand{\epf}{\end{pf}}
\newcommand{\ba}{\begin{array}}
\newcommand{\ea}{\end{array}}
\newcommand{\beq}{\begin{eqnarray}}
\newcommand{\beqq}{\begin{eqnarray*}}
\newcommand{\eeq}{\end{eqnarray}}
\newcommand{\eeqq}{\end{eqnarray*}}

\newcommand{\ds}{\displaystyle}

\newcounter{minutes}
\divide\time by 60
\newcounter{hours}
\multiply\time by 60 \addtocounter{minutes}{-\time}

\begin{document}

\bibliographystyle{amsplain}
\title []
{On certain quasiconformal and elliptic mappings}

\def\thefootnote{}
\footnotetext{ \texttt{\tiny File:~\jobname .tex,
          printed: \number\day-\number\month-\number\year,
          \thehours.\ifnum\theminutes<10{0}\fi\theminutes}
} \makeatletter\def\thefootnote{\@arabic\c@footnote}\makeatother

\author{Shaolin Chen}
 \address{Sh. Chen, College of Mathematics and
Statistics, Hengyang Normal University, Hengyang, Hunan 421008,
People's Republic of China.} \email{mathechen@126.com}

\author{Saminathan Ponnusamy}
\address{S. Ponnusamy, Department of Mathematics,
Indian Institute of Technology Madras, Chennai-600 036, India. }
\email{samy@iitm.ac.in}



\begin{abstract}
Let $\overline{\mathbb{D}}$ be the closure of the unit disk
$\mathbb{D}$ in the complex plane $\mathbb{C}$ and $g$ be a
continuous function in $\overline{\mathbb{D}}$. In this paper, we
discuss some characterizations of elliptic mappings  $f$ satisfying
the Poisson's equation $\Delta f=g$ in $\mathbb{D}$, and then
establish some sharp distortion theorems on elliptic mappings with
the finite perimeter and the finite radial length, respectively. The
obtained results are the extension of the corresponding classical
results.
\end{abstract}

\maketitle \pagestyle{myheadings} \markboth{ Sh. Chen, S.
Ponnusamy}{On certain quasiconformal and elliptic mappings}

\section{Preliminaries and  main results }\label{csw-sec1}

Let  $\mathbb{D}=\{z:\, |z|<1\}$ denote the
open unit disk in the complex plane $\mathbb{C}$ and let $\mathbb{T}=\partial\mathbb{D}$
be the unit circle.
Furthermore, we denote by $\mathcal{C}^{m}(\Omega)$ the set of all complex-valued $m$-times
continuously differentiable functions from $\Omega$ into $\mathbb{C}$, where $\Omega$ is a
subset of $\mathbb{C}$ and $m\in\{0,1,2,\ldots\}$. In particular,
$\mathcal{C}(\Omega):=\mathcal{C}^{0}(\Omega)$ denotes the set of
all continuous functions in $\Omega$.  Let $G$ be a domain of
$\mathbb{C}$ with  $\overline{G}$ be its closure.  We use
$d_{G}(z)$ to denote the Euclidean distance from $z$ to the boundary
$\partial G$ of $G$. Especially, we always use
$d(z)$ for $d_{\mathbb{D}}(z)$.

For a real $2\times2$ matrix $A$,  the matrix norm and the matrix function are defined by
$$\|A\|=\sup\{|Az|:\,|z|=1\}, ~\mbox{ and }~ l(A)=\inf\{|Az|:\,|z|=1\},
$$
respectively. For $z=x+iy\in\mathbb{C}$, the
formal derivative of a complex-valued function $f=u+iv$ is given
by
$$D_{f}:=\left(\begin{array}{cccc}
\ds u_{x}\;~~ u_{y}\\[2mm]
\ds v_{x}\;~~ v_{y}
\end{array}\right),
$$
so that
$$\|D_{f}\|=|f_{z}|+|f_{\overline{z}}| ~\mbox{ and }~ l(D_{f})=\big| |f_{z}|-|f_{\overline{z}}|\big |,
$$
where
\be\label{eq1}
f_{z}=\big( f_x-if_y\big)/2,\;\;\mbox{and}\;\;
f_{\overline{z}}=\big(f_x+if_y\big)/2.
\ee
We use $J_{f}:=\det D_{f} =|f_{z}|^{2}-|f_{\overline{z}}|^{2}$ to denote the {\it Jacobian} of $f$.

A sense-preserving homeomorphism $f$ from a domain $\Omega$ onto $\Omega'$, contained in the {\it Sobolev class}
$W_{loc}^{1,2}(\Omega)$, is said to be a {\it $K$-quasiconformal
mapping} if, for $z\in\Omega$,
$$\|D_{f}(z)\|^{2}\leq K\big | \det D_{f}(z)\big |,~\mbox{i.e.,}~\|D_{f}(z)\|\leq Kl(D_{f}(z)),
$$
where $K\geq1$.

A  mapping $f\in\mathcal{C}^{1}(\Omega)$ is said to be an {\it
elliptic mapping} (or {\it $(K,K')$-elliptic mapping}) if there are
constants $K\geq1$ and $K'\geq0$ such that $f$ satisfies the
following partial differential inequality
$$\|D_{f}(z)\|^{2}\leq KJ_{f}(z)+K'
$$
in the domain $\Omega\subset\mathbb{C}$ (see \cite{FS,Ni}).
Let  $\Omega_{1}$ and $\Omega_{2}$ be two subdomains of $\mathbb{C}$. A
sense-preserving homeomorphism $f:\,\Omega_{1}\rightarrow\Omega_{2}$
is said to be a {\it $(K,K')$-quasiconformal mapping} if $f$ is absolutely continuous
on lines in $\Omega_{1}$, and there are constants $K\geq1$ and $K'\geq0$ such that
$$\|D_{f}(z)\|^{2}\leq KJ_{f}(z)+K',~z\in\Omega_{1}.
$$
Obviously, a $(K,K')$-quasiconformal mapping $f\in\mathcal{C}^{1}(\Omega)$
is an  elliptic mapping in $\Omega$. In particular, if $K'=0$, then the $(K,K')$-quasiconformal
mappings are  $K$-quasiconformal  (cf. \cite{CLSW,K8}).  Moreover, if a
$(K,K')$-quasiconformal mapping is harmonic, then it is said to be
harmonic $(K,K')$-quasiconformal.

In 1968,  Martio \cite{Mar} has discussed the conditions for the
$K$-quasiconformality of  harmonic mappings from the closed unit
disk onto itself.   The harmonic $K$-quasiconformal mappings of
Riemannian manifolds have been considered by Goldberg and Ishihara
(see \cite{GI-1,GI-2}). Tam and Wan \cite{TW} have investigated some
properties of harmonic quasiconformal diffeomorphism and the
universal Teichmuller space. In 2002, Pavlovi\'c \cite{Pav-1} has
generalized the corresponding results of Martio \cite{Mar}. Kalaj
\cite{K-2003}, Partyka and Sakan \cite{PS-2014} have investigated
the $K$-quasiconformality of  harmonic mappings from the  unit disk
onto bounded convex domains. See \cite{K1,K-2015,M-2017,PS,Zhu} and
the references therein for detailed discussions on this topic.
 Recent
papers \cite{CLSW}, \cite{K8} and \cite{K-2019} bring much attention
on the topic of $(K,K')$-quasiconformal mappings in the plane.  This
paper continues the study of previous work of
\cite{CPR,CPR-Nolinear} and is mainly motivated by the articles of
Finn and Serrin \cite{FS}, and Kalaj and Mateljevi\'c \cite{K8}. In
order to state our main results, we need to recall some basic
definitions and some results which motivate the present work.

For $\theta\in[0,2\pi]$ and $z, w\in\mathbb{D}$ with $z\neq w$, let
$$G(z,w)=\log\left|\frac{1-z\overline{w}}{z-w}\right|~\mbox{and}~P(z,e^{i\theta})
=\frac{1-|z|^{2}}{|1-ze^{-i\theta}|^{2}}
$$
denote the {\it Green function} and  the {\it (harmonic) Poisson kernel}, respectively.

Let $\psi: \,\mathbb{T}\rightarrow\mathbb{C}$ be a bounded integrable
function,  and $g\in\mathcal{C}(\overline{\mathbb{D}})$. The
solution to the {\it Poisson equation}
$$ \begin{cases}
\displaystyle \Delta f=g &\, \mbox{in}~\mathbb{D},\\
\displaystyle f=\psi\in L^{1}(\mathbb{T})  &\,
\mbox{in}~\mathbb{T},
\end{cases}
$$
is given by
\beqq\label{eq-1.0}
f(z)=P[\psi](z)-G[g](z),
\eeqq
where
\be\label{eq-2.0}
G[g](z)=\frac{1}{2\pi}\int_{\mathbb{D}}G(z,w)g(w)dA(w),~
~P[\psi](z)=\frac{1}{2\pi}\int_{0}^{2\pi}P(z,e^{it})\psi(e^{it})dt,
\ee
and $dA(w)$ denotes the Lebesgue measure on $\mathbb{D}$. It is well
known that if $\psi$ and $g$ are continuous in $\mathbb{T}$ and in
$\overline{\mathbb{D}}$, respectively, then $f=P[\psi]-G[g]$ has a
continuous extension $\tilde{f}$ to the boundary, and
$\tilde{f}=\psi$ in $\mathbb{T}$ (see \cite[pp. 118-120]{Ho}  and
\cite{K1,K2}).

A continuous increasing function $\omega:\, [0,\infty)\rightarrow
[0,\infty)$ with $\omega(0)=0$ is called a {\it majorant} if
$\omega(t)/t$ is non-increasing for $t>0$ (see
\cite{Dy1,Dy2,P-1999}). For $\alpha>0$ and a majorant $\omega$, we
use $\mathcal{B}_{\omega}^{\alpha}(\mathbb{D})$ to denote the
generalized Bloch-type space of all functions
$f\in\mathcal{C}^{1}(\mathbb{D})$ with
$\|f\|_{\mathcal{B}_{\omega}^{\alpha}(\mathbb{D})}<\infty,$ where
$$\|f\|_{\mathcal{B}_{\omega}^{\alpha}(\mathbb{D})}
=|f(0)|+\sup_{z\in\mathbb{D}}\left\{\|D_{f}(z)\|\omega(d^{\alpha}(z))\right\}.
$$

For a given bounded integrable function $\psi\in L^{1}(\mathbb{T})$
and a given $g\in\mathcal{C}(\overline{\mathbb{D}})$, let
$$\mathcal{F}_{g,\psi}(\mathbb{D})=\{f\in\mathcal{C}^{2}(\mathbb{D}): \, \Delta f=g~\mbox{in}~\mathbb{D}~\mbox{and}~
f=\psi~\mbox{in}~\mathbb{T} \}.
$$
Clearly, functions in $\mathcal{F}_{0,\psi}(\mathbb{D})$ are
harmonic in $\mathbb{D}$. Furthermore, if
$f\in\mathcal{F}_{g,\psi}(\mathbb{D})$, then  $f+G[g]$ is harmonic
in $\mathbb{D}$, and thus, has the representation \be\label{Du}
f+G[g]=h_{1}+\overline{h}_{2}, \ee where  $h_{1}$ and $h_{2}$ are
analytic in $\mathbb{D}$ (cf. \cite{Du}), and $G[g]$ is defined in
{\rm (\ref{eq-2.0})}.

In \cite{K8}, Kalaj and  Mateljevi\'c proved that a harmonic
diffeomorphism between two bounded Jordan domains with
$\mathcal{C}^{2}$ boundaries is a harmonic $(K,K')$-quasiconformal
mapping for some constants $K\geq1$ and $K'\geq0$ if and only if it
is Lipschitz continuous (see \cite[Theorem 1.1 and Corollary
1.3]{K8}). This result can be considered as an extension of the
corresponding results of Martio \cite{Mar}, Pavlovi\'c \cite{Pav-1},
and Partyka and  Sakan \cite{PS}. For related investigations on
this topic, we refer to \cite{PS-2014,Zhu}. In the following, we
will give some characterizations of elliptic mappings without the
$\mathcal{C}^{2}$ boundary hypothesis of the image domains.


\begin{thm}\label{thm-2}
Suppose that $\omega$ is a given majorant.
For a given bounded integrable function $\psi\in L^{1}(\mathbb{T})$
and a given $g\in\mathcal{C}(\overline{\mathbb{D}})$, let
$f\in\mathcal{F}_{g,\psi}(\mathbb{D})$ be a univalent and
sense-preserving  mapping, and $f(\mathbb{D})$ be a convex domain. If there
are two positive constants $C_{1}$,  $C_{2}$ and $\alpha\in[0,1]$ such that for
any $z_{1}, z_{2}\in\mathbb{D}$ with $z_{1}\neq z_{2}$,
\be\label{bi}
\frac{\omega\left(\big((1+|z_{1}|)(1+|z_{2}|)\big)^{\frac{1-\alpha}{2}}\right)}{C_{1}}
\leq\frac{|f(z_{1})-f(z_{2})|}{|z_{1}-z_{2}|}
\leq\frac{C_{1}}{\omega\left(\big(d(z_{1})d(z_{2})\big)^{\frac{1-\alpha}{2}}\right)}
\ee
and
$$\int_{0}^{1}\frac{dt}{\omega\left((d(\Phi(t)))^{1-\alpha}\right)}\leq C_{2},
$$
then $f$ is an elliptic mapping, where $\Phi(t):=f^{-1}(f(z_{1})+t(f(z_{2})-f(z_{1})))$.
\end{thm}

\begin{rem}
In particular, if $\alpha=1$ in (\ref{bi}), then  $f$ is
bi-Lipschitz. In this situation, Theorem \ref{thm-2} is trivial
because all univalent and sense-preserving bi-Lipschitz mappings are
quasiconformal mappings (see Chapter 14.78 in \cite{HKM}). However,
quasiconformal mappings are not necessarily bi-Lipschitz, not even
Lipschitz (see Example \ref{CKW-0}).
\end{rem}

\begin{exam}\label{CKW-0}
Let 
$$ f(z)=\begin{cases}
\displaystyle z\log^{\alpha}\left(\frac{e}{|z|^{2}}\right) & \mbox{ for  $z\in\mathbb{D}\setminus\{0\}$},\\
\displaystyle 0  & \mbox{ for $z=0$},
\end{cases}
$$
where $\alpha\in(0,1/2)$ is a constant. Then $f$ is a quasiconformal
self-homeomorphism of $\mathbb{D}$. However, $f$ is not Lipschitz at
the origin (cf. \cite{KS}).
\end{exam}

\begin{prop}\label{p-1}
For a given bounded integrable function $\psi\in L^{1}(\mathbb{T})$
and a given $g\in\mathcal{C}(\overline{\mathbb{D}})$, let
$f=h_{1}+\overline{h}_{2}-G[g]\in\mathcal{F}_{g,\psi}(\mathbb{D})$,
where $h_{1}$ and $h_{2}$ are analytic in
$\mathbb{D}$. 
If there is a constant $C_{3}\in[1,2)$ such that for any $z_{1},$
$z_{2}\in\mathbb{D}$,
\be\label{c-po-1}
|f(z_{1})-f(z_{2})|\leq C_{3}|h_{1}(z_{1})-h_{1}(z_{2})|,
\ee
then $f$ is an elliptic mapping.
\end{prop}

We remark that the  inverse of Proposition \ref{p-1} does not
necessarily hold (see Example \ref{cexamp-2}).


\begin{exam}\label{cexamp-2}
For $z\in\mathbb{D}$, let $f(z)=3z|z|^{2}-z|z|^{8}$. Then

\begin{enumerate}
\item[{\rm (1)}] $f$ is a $(1,729/(2^{16/3}))$-quasiconformal mapping in $\mathbb{D}$;

\item[{\rm (2)}] $f$ is not a $K$-quasiconformal mapping for any
$K\geq1$;

\item[{\rm (3)}] $f$ does not satisfy the inequality (\ref{c-po-1}).
\end{enumerate}

\end{exam}
\bpf We first prove the univalence of $f$. Suppose on the contrary
that $f$ is not univalent. Then there are two distinct points
$z_{1},~z_{2}\in\mathbb{D}$ such that $f(z_{1})=f(z_{2}),$  which
implies that
\be\label{examp-1}
z_{1}|z_{1}|^{2}(3-|z_{1}|^{6})=z_{2}|z_{2}|^{2}(3-|z_{2}|^{6}).
\ee

$\mathbf{Case ~1.}$ If $|z_{1}|=|z_{2}|$, then, by (\ref{examp-1}),
we have $z_{1}=z_{2}$. This is a contradiction with the assumption.

$\mathbf{Case ~2.}$ If $|z_{1}|\neq|z_{2}|$, then (\ref{examp-1})
reduces to $|z_{1}|^{3}(3-|z_{1}|^{6})=|z_{2}|^{3}(3-|z_{2}|^{6}),$
and consequently
$$(|z_{1}|^{3}-|z_{2}|^{3})(3-|z_{1}|^{6}-|z_{2}|^{6}-|z_{1}|^{3}|z_{2}|^{3})=0.
$$
This implies $|z_{1}|=|z_{2}|=1$, which violates the hypothesis.
Hence $f$ is univalent.

Also, for $z\in\mathbb{D}$, elementary calculations lead to
$$f_{z}(z)=|z|^{2}(6-5|z|^{6})~\mbox{and}~f_{\overline{z}}(z)=z^{2}(3-4|z|^{6}),
$$
which give that
\be\label{examp-2}
\|D_{f}(z)\|^{2}-J_{f}(z)\leq\|D_{f}(z)\|^{2}=\max_{|z|\in[0,1)}\left\{81|z|^{4}(1-|z|^{6})^{2}\right\}=729/(2^{16/3})
\ee
and
\beq\label{examp-3}
\lim_{|z|\rightarrow1^{-}}\frac{|f_{\overline{z}}(z)|}{|f_{z}(z)|}=1.
\eeq
The first assertion and the second assertion easily follows from (\ref{examp-2}) and (\ref{examp-3}), respectively.
The last assertion is obvious.
\epf

For $r\in (0,1)$ and $f\in\mathcal{C}^{1}(\mathbb{D})$,   the {\it perimeter} of the curve
$C(r)=\big\{w=f(re^{i\theta}):\, \theta\in[0,2\pi]\big\}$, with counting
multiplicity, is defined by (cf. \cite{CLP-2017,CPR-Nolinear,C-2019})
\beqq\label{eq-c-1}
\ell_{f}(r)=\int_{0}^{2\pi}|df(re^{i\theta})|=
r\int_{0}^{2\pi}\left|f_{z}(re^{i\theta})-e^{-2i\theta}
f_{\overline{z}}(re^{i\theta})\right|d\theta .
\eeqq
In particular, let $\ell_{f}(1)=\sup_{0<r<1}\ell_{f}(r)$.
Let us recall the following distortion theorem for $K$-quasiconformal harmonic mappings with finite
perimeter.

\begin{Thm}
{\rm (\cite[Theorem 2]{CPR-Nolinear})}\label{CPR-N2015}
Let
$f(z)=\sum_{n=0}^{\infty}a_{n}z^{n}+\sum_{n=1}^{\infty}\overline{b}_{n}\overline{z}^{n}$
be a $K$-quasiconformal harmonic mapping. If $\ell_{f}(1)<\infty$,
then for $n\geq1$, \be\label{CRP-1c}
|a_{n}|+|b_{n}|\leq\frac{K\ell_{f}(1)}{2n\pi},\ee
\be\label{CRP-2c}
\|D_{f}(z)\|\leq\frac{\ell_{f}(1)\sqrt{K}}{2\pi(1-|z|)} \ee and
$f\in\mathcal{B}_{\omega}^{1}(\mathbb{D})$, where $\omega(t)=t$. In
particular, if $K=1,$ then the above two estimates are sharp, and
the extreme function is $f(z)=z.$
\end{Thm}

Concerning the generalized form of inequality (\ref{CRP-1c}),
Mateljevi\'c \cite{MM-2019} proved the following result: Let
$f(z)=\sum_{n=0}^{\infty}a_{n}z^{n}+\sum_{n=1}^{\infty}\overline{b}_{n}\overline{z}^{n}$
be a  harmonic mapping with $\ell_{f}(1)<\infty$. Then, for
$n\geq1$, the inequality
\be\label{Mat-1}|a_{n}|+|b_{n}|\leq\frac{\ell_{f}(1)}{n\pi}\ee holds
 (see \cite[Theorem 10]{MM-2019}).
 Moreover, Kalaj \cite{Kaj-2019}
improved the inequality (\ref{CRP-2c}) and obtained a sharp
inequality for harmonic diffeomorphisms of $\mathbb{D}$. It reads as
follows.


\begin{Thm}\label{Kalaj-2019} If $f$ is a harmonic sense-preserving diffeomorphism of $\mathbb{D}$ onto a Jordan domain $\Omega$ with rectifiable boundary of length $2\pi R$,
then the sharp inequality
\be\label{kalaj-1} |f_{z}(z)|\leq\frac{R}{1-|z|^{2}},~z\in\mathbb{D}
\ee holds, where $R$ is a positive constant. If the equality in
\eqref{kalaj-1} is attained for some $a$, then $\Omega$ is convex,
and there is a holomorphic function
$\mu:~\mathbb{D}\rightarrow\mathbb{D}$ and a constant
$\theta\in[0,2\pi]$, such that \be\label{kalaj-2}
F(z):=e^{-i\theta}f\left(\frac{z+a}{1+z\overline{a}}\right)=R\left(\int_{0}^{z}\frac{dt}{1+t^{2}\mu(t)}+
\overline{\int_{0}^{z}\frac{\mu(t)dt}{1+t^{2}\mu(t)}}\right).\ee
Moreover, every function $f$ defined by {\rm (\ref{kalaj-2})} is a
harmonic diffeomorphism and maps $\mathbb{D}$ to a Jordan domain
bounded by a convex curve of length $2\pi R$ and the inequality
\eqref{kalaj-1} is attained for $z=a$.
\end{Thm}


The following result is also a generalization of Theorem
\Ref{CPR-N2015}.



\begin{thm}\label{thm-3} Let $K\geq1,$  $K'\geq0$ and $R>0$ be
constants. If $f(z)=\sum_{n=0}^{\infty}a_{n}z^{n}+\sum_{n=1}^{\infty}\overline{b}_{n}\overline{z}^{n}$ is a harmonic $(K,K')$-quasiconformal mapping of
$\mathbb{D}$ onto a Jordan domain $\Omega$ with rectifiable boundary
of length $2\pi R$, then for $n\geq1$, \be\label{chen-1.0} |a_{n}|+|b_{n}|\leq
\frac{\sqrt{K'}+KR}{n}, \ee
\be\label{CP-K}
\|D_{f}(z)\|\leq\left(R+\frac{-R+\sqrt{K'+KK'+K^{2}R^{2}}}{1+K}\right)\frac{1}{1-|z|^{2}}, ~\mbox{ for $z\in\mathbb{D}$},
\ee 
and $f\in\mathcal{B}_{\omega}^{1}(\mathbb{D})$.

 In particular, if $K'=K-1=0$ and $\Omega=\mathbb{D},$ then the
estimates of {\rm (\ref{chen-1.0})} is
sharp, and the extreme function is $f(z)=z$ for $z\in\mathbb{D}$.
Moreover, if $K'=K-1=0$ and $\Omega=\mathbb{D},$ then the equal sign
occurs in {\rm (\ref{CP-K})} for some fixed $z=a$ if and only if
$f(z)=e^{it}\frac{z-a}{1-z\overline{a}}$, where $t\in[0,2\pi].$
\end{thm}
We remark that if $K'=0$ and $K\in[1,2]$, then the
inequality (\ref{chen-1.0}) is better than (\ref{Mat-1}).

Let $f\in\mathcal{C}^{1}(\mathbb{D})$. Then, for
$\theta\in[0,2\pi]$, the {\it radial length} of the curve
$C_{\theta}(r)=\big\{w=f(\rho e^{i\theta}):\, 0\leq\rho\leq
r<1\big\}$, with counting multiplicity, is defined by

\beqq\label{eq-cp-1}
\ell_{f}^{\ast}(r,\theta)
=\int_{0}^{r}|df(\rho e^{i\theta})|
= \int_{0}^{r}\left|f_{z}(\rho e^{i\theta})+e^{-2i\theta} f_{\overline{z}}(\rho
e^{i\theta})\right|d\rho.
\eeqq
In particular, let
$$\ell_{f}^{\ast}(1,\theta)=\sup_{0\leq r<1}\ell_{f}^{\ast}(r,\theta).
$$
We refer the reader to \cite{CPR-Nolinear, CP-2019} for some
discussion of the radial length. In particular, the following result
establishes the Fourier coefficient estimates of $K$-quasiconformal
harmonic mappings with the finite radial length.

\begin{Thm}{\rm (\cite[Theorem 4]{CLP-2017})}\label{CLP-2017}
Let $f(z)=\sum_{n=0}^{\infty}a_{n}z^{n}+\sum_{n=1}^{\infty}\overline{b}_{n}\overline{z}^{n}$
be a harmonic $K$-quasiconformal  mapping in $\mathbb{D}$. If $\ell_{f}^{\ast}(1)=\sup_{\theta\in[0,2\pi]}\ell_{f}^{\ast}(1,\theta)<\infty$,
then
\be\label{eq-2017}
|a_{n}|+|b_{n}|\leq K\ell_{f}^{\ast}(1)~\mbox{for}~n\geq1.
\ee
Moreover, if $K=1$, then the estimate {\rm (\ref{eq-2017})} is sharp and the extreme function
is $f(z)=\ell_{f}^{\ast}(1)z$ for $z\in\mathbb{D}$.
\end{Thm}

We improve Theorem \Ref{CLP-2017} into the following form.

\begin{thm}\label{thm-4}
For $K\geq1$ and $K'\geq0$, let $f(z)=\sum_{n=0}^{\infty}a_{n}z^{n}+\sum_{n=1}^{\infty}\overline{b}_{n}\overline{z}^{n}$
be a harmonic $(K,K')$-quasiconformal mapping in $\mathbb{D}$. If
$\ell_{f}^{\ast}(1)=\sup_{\theta\in[0,2\pi]}\ell_{f}^{\ast}(1,\theta)<\infty$,
then for $n\geq1$,

\be\label{chen-1.2}
|a_{n}|+|b_{n}|\leq\sqrt{K'}+K\ell_{f}^{\ast}(1).
\ee
In particular, if $K'=K-1=0$, then the estimate {\rm(\ref{chen-1.2})}
is sharp, and the extreme function is $f(z)=\ell_{f}^{\ast}(1)z$ for $z\in\mathbb{D}$.

\end{thm}

 The proof of Theorems \ref{thm-2}, \ref{thm-3}
and \ref{thm-4}, and Proposition \ref{p-1} will be presented in
Section \ref{csw-sec2}.






\section{The proofs of the  main results }\label{csw-sec2}

In this section, we shall prove Theorems \ref{thm-2}
and \ref{thm-4}, and Propositions \ref{p-1} and \ref{thm-3}. We start with some
useful Lemmas.

\begin{Lem}{\rm (\cite[Lemma 6]{CPR})}\label{Lem-A}
Let $\omega$ be a majorant and $\nu\in[0,1]$. Then for
$t\in[0,\infty]$, $\omega(\nu t)\geq \nu\omega(t)$.
\end{Lem}

\begin{lem}\label{thm-1}
Let $\omega$ be a  majorant and $\alpha\in[0,1]$ be a constant. Suppose that
$f\in\mathcal{C}^{1}(\mathbb{D})$ is univalent, and $f(\mathbb{D})$
is a convex domain. Then the following two statements are
equivalent:
\begin{enumerate}
\item[{\rm (a)}] For any $z_{1},$ $z_{2}$ with $z_{1}\neq z_{2}$,
there exists a constant $C_{4}$ such that
$$\frac{1}{C_{4}}\omega\left(\big((1+|z_{1}|)(1+|z_{2}|)\big)^{(1-\alpha)/2}\right)
\leq\frac{|f(z_{1})-f(z_{2})|}{|z_{1}-z_{2}|}
\leq\frac{C_{4}}{\omega\left(\big(d(z_{1})d(z_{2})\big)^{(1-\alpha)/2}\right)}.
$$

\item[{\rm (b)}] For any $z\in\mathbb{D}$, there is a constant
$C_{5}>0$ such that
$$\frac{1}{C_{5}}\omega\left((1+|z|)^{1-\alpha}\right)\leq l(D_{f}(z))
\leq\|D_{f}(z)\|\leq\frac{C_{5}}{\omega\left((d(z))^{1-\alpha}\right)}.
$$
\end{enumerate}
\end{lem}
\bpf We first prove (a)$~\Rightarrow~$(b). For $\theta\in[0,2\pi]$ and
$z=x+iy\in\mathbb{D}$, elementary computations lead to (see \eqref{eq1})
$$ f_{z}(z) + f_{\overline{z}}(z) =f_x(z)\;\;\mbox{and}\;\;
i(f_z(z)-f_{\overline{z}}(z))=f_y(z)
$$
and therefore,
\beq\label{e-1}
f_{x}(z)\cos\theta+f_{y}(z)\sin\theta
&=&f_{z}(z)e^{i\theta}+f_{\overline{z}}(z)e^{-i\theta}.
\eeq
For $r\in[0,1-|z|)$, let $\xi=z+re^{i\theta}$. Then, by (a) and (\ref{e-1}), we obtain
\beqq
\|D_{f}(z)\|&=&\max_{\theta\in[0,2\pi]}\left|f_{x}(z)\cos\theta+f_{y}(z)\sin\theta\right|
=\max_{\theta\in[0,2\pi]}\left\{\lim_{r\rightarrow0^{+}}\frac{|f(z)-f(\xi)|}{|z-\xi|}\right\}\\
&\leq&\max_{\theta\in[0,2\pi]}\left\{\lim_{r\rightarrow0^{+}}\frac{C_{4}}
{\omega\left(\big(d(z)d(z+re^{i\theta})\big)^{\frac{1-\alpha}{2}}\right)}\right\}
=\frac{C_{4}}{\omega\left((d(z))^{1-\alpha}\right)}
\eeqq
and
\beqq
\|D_{f}(z)\|&\geq&l(D_{f}(z))=\min_{\theta\in[0,2\pi]}\left|f_{x}(z)\cos\theta+f_{y}(z)\sin\theta\right|\\
&=&\min_{\theta\in[0,2\pi]}\left\{\lim_{r\rightarrow0^{+}}\frac{|f(z)-f(\xi)|}{|z-\xi|}\right\}\\
&\geq&\min_{\theta\in[0,2\pi]}\left\{\lim_{r\rightarrow0^{+}}\frac{\omega\left(\big((1+|z|)(1+|z+re^{i\theta}|)\big)^{\frac{1-\alpha}{2}}\right)}{C_{4}}\right\}\\
&=&\frac{\omega\left((1+|z|)^{1-\alpha}\right)}{C_{4}}.
\eeqq

Now we prove (b)$~\Rightarrow~$(a). For $t\in[0,1]$ and $z_{1},
z_{2}\in\mathbb{D}$ with $z_{1}\neq z_{2}$, let
$$\chi(t)=z_{1}t+(1-t)z_{2}.
$$
Since
$$d(\chi(t))\geq1-t|z_{1}|-(1-t)|z_{2}|\geq1-t-(1-t)|z_{2}|=(1-t)d(z_{2})
$$
and
$$d(\chi(t))\geq1-t|z_{1}|-(1-t)|z_{2}|\geq1-t|z_{1}|-(1-t)=td(z_{1}),
$$
we see that
\be\label{ch-1.1}\big(d(\chi(t))\big)^{1-\alpha}\geq\big(t(1-t)d(z_{1})d(z_{2})\big)^{\frac{1-\alpha}{2}}.
\ee
By calculations, we have
\beq\label{ch-1.2}
|f(z_{1})-f(z_{2})|&=&\left|\int_{0}^{1}\big(f_{w}(\chi(t))\chi'(t)+f_{\overline{w}}(\chi(t))\overline{\chi'(t)}\big)dt\right|\\
\nonumber &\leq&|z_{1}-z_{2}|\int_{0}^{1}\|D_{f}(\chi(t))\|dt,
\eeq
where $w=\chi(t)$.

Now, we estimate the integral on the right. By (b) and (\ref{ch-1.1}), we have
\beqq
\int_{0}^{1}\|D_{f}(\chi(t))\|dt\leq\int_{0}^{1}\frac{C_{5}dt}{\omega\left(\big(d(\chi(t))\big)^{1-\alpha}\right)}\leq
\int_{0}^{1}\frac{C_{5}dt}{\omega\left(\big(t(1-t)d(z_{1})d(z_{2})\big)^{\frac{1-\alpha}{2}}\right)},
\eeqq
which, together with Lemma \Ref{Lem-A}, implies that
\beq\label{ch-1.3}
\int_{0}^{1}\|D_{f}(\chi(t))\|dt
&\leq&\frac{C_{5}}{\omega\left(\big(d(z_{1})d(z_{2})\big)^{\frac{1-\alpha}{2}}\right)}\int_{0}^{1}\frac{dt}{t^{\frac{1-\alpha}{2}}(1-t)^{\frac{1-\alpha}{2}}}\\
\nonumber
 &=&\frac{\Gamma^{2}\left(\frac{1+\alpha}{2}\right)}{\Gamma(1+\alpha)}\frac{C_{5}}{\omega\left(\big(d(z_{1})d(z_{2})\big)^{\frac{1-\alpha}{2}}\right)},
\eeq
where $\Gamma$ denotes the usual Gamma function.

It follows from (\ref{ch-1.2}) and (\ref{ch-1.3}) that
$$\frac{|f(z_{1})-f(z_{2})|}{|z_{1}-z_{2}|}
\leq\frac{\Gamma^{2}\left(\frac{1+\alpha}{2}\right)}{\Gamma(1+\alpha)}\frac{C_{5}}{\omega\left(\big(d(z_{1})d(z_{2})\big)^{\frac{1-\alpha}{2}}\right)}.
$$

On the other hand, for any $z_{1}, z_{2}\in\mathbb{D}$ with $z_{1}\neq z_{2}$, let $w_{1}=f(z_{1})$ and $w_{2}=f(z_{2})$. For
$t\in[0,1]$, let
\be\label{ch-1.4}
\gamma(t)=tw_{1}+(1-t)w_{2}
\ee
be the straight line segment connecting $w_{1}$ and $w_{2}$. Since $f(\mathbb{D})$
is a convex domain, we see that $\gamma(t)\subset f(\mathbb{D})$ and
$\eta(t)=f^{-1}(\gamma(t))\subset\mathbb{D}$ for $t\in[0,1]$. It is
not difficult to know that
$$1+|\eta(t)|\geq\frac{1+|z_{1}|}{2}~\mbox{and}~1+|\eta(t)|\geq\frac{1+|z_{2}|}{2},
$$
which, together with Lemma \Ref{Lem-A}, implies that
\beq\label{ch-1.5}
\omega\left((1+|\eta(t)|)^{1-\alpha}\right)&\geq&\omega\left(\frac{\big((1+|z_{1}|)(1+|z_{2}|)\big)^{\frac{1-\alpha}{2}}}{2^{1-\alpha}}\right)\\
\nonumber
&\geq&\frac{\omega\left(\big((1+|z_{1}|)(1+|z_{2}|)\big)^{\frac{1-\alpha}{2}}\right)}{2^{1-\alpha}}.
\eeq
By (b) and (\ref{ch-1.5}), we obtain
\beq\label{ch-1.6}
\int_{0}^{1}l(D_{f}(\eta(t)))|\eta'(t)|dt&\geq&\frac{1}{C_{5}}\int_{0}^{1}\omega\left((1+|\eta(t)|)^{1-\alpha}\right)|\eta'(t)|dt\\
\nonumber
&\geq&\frac{\omega\left(\big((1+|z_{1}|)(1+|z_{2}|)\big)^{\frac{1-\alpha}{2}}\right)}{2^{1-\alpha}C_{5}}\int_{0}^{1}|\eta'(t)|dt\\
\nonumber
&\geq&\frac{\omega\left(\big((1+|z_{1}|)(1+|z_{2}|)\big)^{\frac{1-\alpha}{2}}\right)}{2^{1-\alpha}C_{5}}|z_{1}-z_{2}|,
\eeq

It follows from (\ref{ch-1.4}) that
\beqq
|w_{1}-w_{2}|&=&|\gamma'(t)|=\int_{0}^{1}|\gamma'(t)|dt\\
\nonumber
&=&\int_{0}^{1}\left|f_{\zeta}(\eta(t))\eta'(t)+f_{\overline{\zeta}}(\eta(t))\overline{\eta'(t)}\right|dt\\
\nonumber &\geq&\int_{0}^{1}l(D_{f}(\eta(t)))|\eta'(t)|dt,
\eeqq
which, together with (\ref{ch-1.6}), yields that
$$\frac{|w_{1}-w_{2}|}{|z_{1}-z_{2}|}
\geq\frac{\omega\left(\big((1+|z_{1}|)(1+|z_{2}|)\big)^{\frac{1-\alpha}{2}}\right)}{2^{1-\alpha}C_{5}}.
$$
The proof of this lemma is finished. \epf


\begin{Lem}{\rm (\cite[Lemma 2.7]{K2})}\label{KP-2011}
If  $g\in\mathcal{C}(\overline{\mathbb{D}})$, then
$$\max\left\{\left|\frac{\partial}{\partial z}G[g](z)\right|,~\left|\frac{\partial}{\partial \overline{z}}G[g](z)\right|\right\}
\leq\frac{1}{3}\|g\|_{\infty}~\mbox{for}~z\in\mathbb{D}.
$$
\end{Lem}

\begin{lem}\label{lem-1.0}
Let $f\in\mathcal{C}^{1}(\mathbb{D})$ be a sense-preserving mapping.
Then $f$ is an elliptic mapping if and only if  there exist
constants $k_{1}\in[0,1)$ and $k_{2}\in[0,\infty)$ such that
\be\label{ch-01}|f_{\overline{z}}(z)|
\leq k_{1}|f_{z}(z)|+k_{2}~\mbox{for}~z\in\mathbb{D}.
\ee
\end{lem}
\bpf We first prove the sufficiency.  By (\ref{ch-01}), for
$z\in\mathbb{D}$, we have
\beq\label{ch-02}
\|D_{f}(z)\|&\leq&\left(\frac{1+k_{1}}{1-k_{1}}\right)l(D_{f}(z))+\frac{2k_{2}}{1-k_{1}}\\
\nonumber
&\leq&\left(\frac{1+k_{1}}{1-k_{1}}\right)l(D_{f}(z))+\sqrt{\left(\frac{1+k_{1}}{1-k_{1}}\right)^{2}l^{2}(D_{f}(z))+\frac{4k_{2}^{2}}{(1-k_{1})^{2}}}.
\eeq

$\mathbf{Case ~1.}$ For all $z\in \mathbb{D}$, $\ds \|D_{f}(z)\|\leq\left(\frac{1+k_{1}}{1-k_{1}}\right)l(D_{f}(z)).$

In this case, it is easy to know that
$$\|D_{f}(z)\|^{2}\leq\left(\frac{1+k_{1}}{1-k_{1}}\right)J_{f}(z),
$$
which implies that $f$ is an elliptic mapping.

$\mathbf{Case ~2.}$ There is a subset $E$ of $\mathbb{D}$ such that  
$$\|D_{f}(z)\|>\left(\frac{1+k_{1}}{1-k_{1}}\right)l(D_{f}(z)) ~\mbox{ for $z\in E$}.
$$ 

In this case,  it follows from (\ref{ch-02}) that 

$$\left[\|D_{f}(z)\|-\left(\frac{1+k_{1}}{1-k_{1}}\right)l(D_{f}(z))\right ]^{2}
\leq\left(\frac{1+k_{1}}{1-k_{1}}\right)^{2}l^{2}(D_{f}(z))+\frac{4k_{2}^{2}}{(1-k_{1})^{2}} ~\mbox{ for $z\in E$},
$$
which implies that
\be\label{ch-j1}
\|D_{f}(z)\|^{2}\leq2\left(\frac{1+k_{1}}{1-k_{1}}\right)J_{f}(z)+\frac{4k_{2}^{2}}{(1-k_{1})^{2}}.
\ee
On the other hand, for $z\in\mathbb{D}\backslash E$, we have
$$\|D_{f}(z)\|^{2}\leq\left(\frac{1+k_{1}}{1-k_{1}}\right)J_{f}(z),
$$
which, together with (\ref{ch-j1}), implies that  $f$  is also an elliptic mapping in $\mathbb{D}$.

Next, we show  the necessity. If $f$ is an elliptic mapping, then
there exist constants $K\geq1$ and $K'\geq0$ such that
$$\|D_{f}(z)\|^{2}\leq KJ_{f}(z)+K' ~\mbox{ for $z\in\mathbb{D}$}.
$$
This gives that
\beqq
\|D_{f}(z)\|\leq\frac{Kl(D_{f}(z))+\sqrt{\big(Kl(D_{f}(z))\big)^{2}+4K'}}{2}
\leq Kl(D_{f}(z))+\sqrt{K'},
\eeqq
and consequently
$$|f_{\overline{z}}(z)|\leq\frac{K-1}{K+1}|f_{z}(z)|+\frac{\sqrt{K'}}{1+K}.
$$
The proof of this lemma is complete.
\epf

\subsection*{The proof of Theorem \ref{thm-2}} Differentiating both sides of the equation $f^{-1}(f(z))=z$ and then
simplifying the resulting relations lead to the formulae
\be\label{ch-1.7}
 (f^{-1})_{w}=\frac{\overline{f_{z}}}{J_{f}} ~\mbox{ and }~
 (f^{-1})_{\overline{w}}=-\frac{f_{\overline{z}}}{J_{f}},
\ee
where $w=f(z)$. Next, for  $z_{1}, z_{2}\in\mathbb{D}$ with $z_{1}\neq z_{2}$, we let
$$\varphi(t)=t(f(z_{1})-f(z_{2}))+f(z_{2}),
$$
where $t\in[0,1]$. Since $f(\mathbb{D})$ is a convex domain, we see that
$\varphi(t)\subset f(\mathbb{D})$ and $\Phi(t):=f^{-1}(\varphi(t))\subseteq\mathbb{D}$ for $t\in[0,1]$.
For $z\in\mathbb{D}$, let $F(z)=f(z)+G[g](z)$. Then $F$ is harmonic
in $\mathbb{D}$, and $F$ can be represented by
$F=h_{1}+\overline{h}_{2}$ in $\mathbb{D}$, where $h_{j}~(j=1,2)$
are analytic in $\mathbb{D}$.

With $z=\Phi(t)$, we have $\Phi(0)=z_2$, $\Phi(1)=z_1$ and thus, by (\ref{ch-1.7}), we obtain
\beqq
h_{1}(z_{1})-h_{1}(z_{2})&=&\int_{z_2}^{z_1}h_{1}'(z)dz=\int_{0}^{1}h'_{1}(\Phi(t))\Phi'(t)dt\\
&=&\int_{0}^{1}h'_{1}(\Phi(t))\left(f^{-1}_{w}(\varphi(t))\varphi'(t)+f^{-1}_{\overline{w}}(\varphi(t))\overline{\varphi'(t)}\right)dt\\
&=&\int_{0}^{1}h'_{1}(\Phi(t))\left(\frac{\overline{f_{z}(\Phi(t))}}{J_{f}(\Phi(t))}\varphi'(t)-\frac{f_{\overline{z}}(\Phi(t))}{J_{f}(\Phi(t))}\overline{\varphi'(t)}\right)dt,
\eeqq
which implies that
\beq\label{ch-1.8}
\left|\frac{h_{1}(z_{1})-h_{1}(z_{2})}{f(z_{1})-f(z_{2})}\right|
&=&\left|\int_{0}^{1}\frac{h'_{1}(\Phi(t))\left(\overline{f_{z}(\Phi(t))}-
f_{\overline{z}}(\Phi(t))\frac{\overline{\varphi'(t)}}{\varphi'(t)}\right)}{J_{f}(\Phi(t))}dt\right|\\
\nonumber
&\leq&\int_{0}^{1}\frac{|h'_{1}(\Phi(t))|\|D_{f}(\Phi(t))\|}{J_{f}(\Phi(t))}dt\\
\nonumber
&\leq&\int_{0}^{1}\frac{l(D_{f}(\Phi(t)))+|G[g]_{z}(\Phi(t))|+|f_{\overline{z}}(\Phi(t))|}{l(D_{f}(\Phi(t)))}dt.
\eeq
It follows from Lemma \ref{thm-1} that there exists a positive
constant $C_{6}$ such that

\beqq
\frac{|G[g]_{z}(\Phi(t))|+|f_{\overline{z}}(\Phi(t))|}{l(D_{f}(\Phi(t)))}&\leq&\frac{|G[g]_{z}(\Phi(t))|+\|D_{f}(\Phi(t))\|}{l(D_{f}(\Phi(t)))}\\
\nonumber
&\leq&\frac{C_{6}|G[g]_{z}(\Phi(t))|+\frac{C_{6}}{\omega\left((d(\Phi(t)))^{1-\alpha}\right)}}{\omega\left((1+|\Phi(t)|)^{1-\alpha}\right)}\\
&\leq&\frac{C_{6}}{\omega(1)}\left(|G[g]_{z}(\Phi(t))|+\frac{1}{\omega\left((d(\Phi(t)))^{1-\alpha}\right)}\right),
\eeqq
which, together with (\ref{ch-1.8}) and Lemma \Ref{KP-2011}, gives that
\beq\label{ch-1.9}
\left|\frac{h_{1}(z_{1})-h_{1}(z_{2})}{f(z_{1})-f(z_{2})}\right|&\leq&1+\frac{C_{6}}{\omega(1)}\bigg(\int_{0}^{1}|G[g]_{z}(\Phi(t))|dt\\
\nonumber &&+
\int_{0}^{1}\frac{dt}{\omega\big((d(\Phi(t)))^{1-\alpha}\big)}\bigg)
\leq\mu,
\eeq
where
$$\mu:=1+\frac{C_{6}\|g\|_{\infty}}{3\omega(1)}+\frac{C_{6}C_{2}}{\omega(1)}\geq1.
$$
For $\theta\in[0,2\pi]$ and any fixed point $z\in\mathbb{D}$, let
$\varsigma=z+re^{i\theta}$, where $r\in[0,1-|z|)$. By (\ref{ch-1.9}), we have
\beqq
\frac{1}{\mu}\lim_{r\rightarrow0^{+}}\left|\frac{h_{1}(\varsigma)-h_{1}(z)}{\varsigma-z}\right|&\leq&\lim_{r\rightarrow0^{+}}
\bigg|\frac{h_{1}(\varsigma)-h_{1}(z)}{\varsigma-z}+\frac{\overline{h_{2}(\varsigma)-h_{2}(z)}}{\varsigma-z}+\frac{G[g](z)-G[g](\varsigma)}{\varsigma-z}\bigg|\\
&=&\big|h_{1}'(z)e^{i\theta}+\overline{h_{2}'(z)}e^{-i\theta}-(G[g]_{z}(z)e^{i\theta}+G[g]_{\overline{z}}(z)e^{-i\theta})\big|
\eeqq
which yields that
\beq\label{ch-2.0}
\frac{1}{\mu}|h_{1}'(z)|&\leq&\min_{\theta\in[0,2\pi]}\big|h_{1}'(z)-G[g]_{z}(z)+(\overline{h'_{2}(z)}-G[g]_{\overline{z}}(z))e^{-2i\theta}\big|\\
\nonumber
&=&\min_{\theta\in[0,2\pi]}\big|f_{z}(z)+f_{\overline{z}}(z)e^{-2i\theta}\big|=l(D_{f}(z)).
\eeq
Since
$$|h_{1}'|\geq|h_{1}'-G[g]_{z}|-|G[g]_{z}|=|f_{z}|-|G[g]_{z}|,
$$
by (\ref{ch-2.0}) and Lemma \Ref{KP-2011}, we see that
\beq\label{ch-2.1}
|f_{\overline{z}}|\leq\frac{(\mu-1)}{\mu}|f_{z}|+\frac{1}{3\mu}\|g\|_{\infty}.
\eeq
It follows from (\ref{ch-2.1}) and Lemma \ref{lem-1.0} that $f$
is an elliptic mapping. The proof of this theorem is complete.
\qed

\subsection*{The proof of Proposition \ref{p-1}} For $\theta\in[0,2\pi]$ and any fixed point $z\in\mathbb{D}$,
let $\varsigma=z+re^{i\theta}$, where $r\in[0,1-|z|)$. Then by (\ref{e-1}) and (\ref{c-po-1}), we have
\beqq
\|D_{f}(z)\|&=&\max_{\theta\in[0,2\pi]}\left|f_{x}(z)\cos\theta+f_{y}(z)\sin\theta\right|=\max_{\theta\in[0,2\pi]}\left\{\lim_{r\rightarrow0^{+}}
\left|\frac{f(\varsigma)-f(z)}{\varsigma-z}\right|\right\}\\
&\leq& C_{3}\lim_{r\rightarrow0^{+}}\left|\frac{h_{1}(\varsigma)-h_{1}(z)}{\varsigma-z}\right|=C_{3}|h_{1}'(z)|\\
&\leq&C_{3}|f_{z}(z)|+C_{3}|G[g]_{z}(z)|,
\eeqq
which, together with Lemma \Ref{KP-2011},  yields that
\be\label{ch-2.6h}
|f_{\overline{z}}(z)|\leq (C_{3}-1)|f_{z}(z)|+\frac{C_{3}}{3}\|g\|_{\infty}.
\ee
It follows from (\ref{ch-2.6h}) and Lemma \ref{lem-1.0} that $f$ is an elliptic
mapping. The proof of this proposition is complete.
\qed



\subsection*{The proof of Theorem \ref{thm-3}} We first prove (\ref{chen-1.0}).
Since
$$\|D_{f}(z)\|^{2}\leq KJ_{f}(z)+K',
$$
we see that
\be\label{bbb-1}
\|D_{f}(z)\|\leq\frac{Kl(D_{f}(z))+\sqrt{\big(Kl(D_{f}(z))\big)^{2}+4K'}}{2}
\leq Kl(D_{f}(z))+\sqrt{K'},
\ee
where $z\in\mathbb{D}$.
It follows from (\ref{bbb-1}) that

\beqq
\ell_{f}(r)&=&
r\int_{0}^{2\pi}\left|f_{z}(re^{i\theta})-e^{-2i\theta}
f_{\overline{z}}(re^{i\theta})\right|d\theta \geq r\int_{0}^{2\pi}l(D_{f}(re^{i\theta}))d\theta\\
&\geq& r\int_{0}^{2\pi}\left(\frac{\|D_{f}(re^{i\theta})\|-\sqrt{K'}}{K}\right)d\theta,
\eeqq
which gives that
\be\label{bbb-2}
r\int_{0}^{2\pi}\|D_{f}(re^{i\theta})\|d\theta\leq2\pi r\sqrt{K'}+K\ell_{f}(r).
\ee


For $n\geq1$ and $r\in(0,1)$, it follows from  Cauchy's integral
formula that,
$$na_{n}=\frac{1}{2\pi i}\int_{|z|=r} \frac{\partial f(z)}{\partial z} \frac{dz}{z^{n}}
~\mbox{ and }~nb_{n}=\frac{1}{2\pi i}\int_{|z|=r} \left(\frac{\partial
\overline{f(z)}}{\partial z}\right)\frac{dz}{z^{n}},$$
which, together with  the inequality (\ref{bbb-2}), implies that
\beqq
n(|a_{n}|+|b_{n}|)&\leq&\frac{1}{2\pi r^{n}}\int_{0}^{2\pi}r\|D_{f}(re^{i\theta})\|d\theta\\
&\leq&\frac{1}{2\pi r^{n}}\left(2\pi r\sqrt{K'}+K\ell_{f}(r)\right)\\
&\leq&\frac{1}{2\pi r^{n}}\left(2\pi \sqrt{K'}+K\ell_{f}(1)\right).
\eeqq 
Consequently, \beq\label{ch-2.2c} |a_{n}|+|b_{n}|\leq
\frac{\sqrt{K'}}{n}+\frac{K\ell_{f}(1)}{2n\pi}=\frac{\sqrt{K'}}{n}+\frac{KR}{n}.
\eeq This proves (\ref{chen-1.0})

Next, we prove  (\ref{CP-K}).  For $z\in\mathbb{D}$, it follows from the inequality
$$\|D_{f}(z)\|^{2}\leq KJ_{f}(z)+K'$$ that
$$|f_{\overline{z}}(z)|\leq\frac{-|f_{z}(z)|+\sqrt{K'+KK'+K^{2}|f_{z}(z)|^{2}}}{1+K},$$
which, together with (\ref{kalaj-1}), yields that
\beqq\|D_{f}(z)\|\leq\left(R+\frac{-R+\sqrt{K'+KK'+K^{2}R^{2}}}{1+K}\right)\frac{1}{1-|z|^{2}}.\eeqq
The proof of the theorem is complete.
\qed

\begin{Thm} {\rm  (\cite[Theorem 2]{B})}\label{Thm-cs}
Let $\phi$ be subharmonic in $\mathbb{D}$. If for all $r\in[0,1)$,
$$\mathcal{A}(r)=\sup_{\theta\in[0,2\pi]}\int_{0}^{r}\phi(\rho e^{i\theta})\,d\rho\leq1,
$$
then $\mathcal{A}(r)\leq r$.
\end{Thm}

\subsection*{The proof of Theorem \ref{thm-4}}
Since $f$ is a univalent $(K,K')$-elliptic mapping,  we see that 
\beqq
\|D_{f}(z)\|^{2}\leq K\|D_{f}(z)\|l(D_{f}(z))+K' ~\mbox{ for $z\in\mathbb{D}$}.
\eeqq
This gives that
\be\label{ch-3.0}
\|D_{f}(z)\|\leq\frac{Kl(D_{f}(z))+\sqrt{\big(Kl(D_{f}(z))\big)^{2}+4K'}}{2}
\leq Kl(D_{f}(z))+\sqrt{K'}.
\ee
It follows from (\ref{ch-3.0})  that, for $\theta\in[0,2\pi]$ and  $r\in(0,1)$,
\begin{eqnarray*}
\ell_{f}^{\ast}(r,\theta) &=&\int_{0}^{r}|f_{z}(\rho
e^{i\theta})+e^{-2i\theta}f_{\overline{z}}(\rho e^{i\theta})|\,d\rho
\geq\int_{0}^{r}l(D_f(\rho e^{i\theta}))\,d\rho\\
&\geq&\frac{1}{K}\int_{0}^{r}\left(\|D_{f}(\rho
e^{i\theta})\|-\sqrt{K'}\right)\,d\rho
\end{eqnarray*}
and consequently \beq\label{ch-3.1} \int_{0}^{r}\|D_{f}(\rho
e^{i\theta})\|\,d\rho\leq\sqrt{K'}r+K\ell_{f}^{\ast}(r,\theta)\leq\sqrt{K'}+K\ell_{f}^{\ast}(1).
\eeq Inequality (\ref{ch-3.1}) and Theorem \Ref{Thm-cs} lead to
\beq\label{ch-3.2} \int_{0}^{r}\|D_{f}(\rho
e^{i\theta})\|\,d\rho\leq\left(\sqrt{K'}+K\ell_{f}^{\ast}(1)\right)r.
\eeq
The Cauchy integral formula  shows that, for $\rho\in(0,1)$,
$$na_{n}=\frac{1}{2\pi i}\int_{|z|=\rho} \frac{\partial f(z)}{\partial z}\frac{dz}{z^{n}}
~\mbox{and}~nb_{n}=\frac{1}{2\pi i}\int_{|z|=\rho}
\left(\frac{\partial \overline{f(z)}}{\partial
z}\right)\frac{dz}{z^{n}},
$$
which yields that \beq\label{ch-3.3} n(|a_{n}|+|b_{n}|)
&\leq&\frac{1}{2\pi \rho^{n}}\int_{0}^{2\pi}\rho\|D_{f}(\rho
e^{i\theta})\|d\theta. \eeq Then combining (\ref{ch-3.2}) and
(\ref{ch-3.3}) gives the final estimate for $|a_{n}|+|b_{n}|$,
namely, \beqq 2\pi
n(|a_{n}|+|b_{n}|)\int_{0}^{r}\rho^{n-1}d\rho&\leq&\int_{0}^{2\pi}\left(\int_{0}^{r}\|D_{f}(\rho
e^{i\theta})\|\,d\rho\right)\,d\theta\\
&\leq&\int_{0}^{2\pi}\left(\sqrt{K'}+K\ell_{f}^{\ast}(1)\right)r\,d\theta
\eeqq and consequently \beqq
|a_{n}|+|b_{n}|\leq\inf_{r\in(0,1)}\left(\frac{\sqrt{K'}+K\ell_{f}^{\ast}(1)}{r^{n-1}}\right)
=\sqrt{K'}+K\ell_{f}^{\ast}(1). \eeqq The proof of this theorem is
complete. \qed







\bigskip

{\bf Acknowledgements:} We thank  anonymous referees for   drawing our attention to
references \cite{Kaj-2019,MM-2019} and providing us with a lot of valuable comments.
This
research was partly supported by the Hunan Provincial Education
Department Outstanding Youth Project (No. 18B365), the Science and
Technology Plan Project of Hengyang City (No. 2018KJ125), the
exchange project for the third regular session of the
China-Montenegro Committee for Cooperation in Science and Technology
(No. 3-13), the
Science and Technology Plan Project of Hunan Province (No.
2016TP1020), the Science and Technology Plan Project of Hengyang
City (No. 2017KJ183),  the Application-Oriented Characterized
Disciplines,  Double First-Class University Project of Hunan
Province (Xiangjiaotong [2018]469), and the Open Fund Project of
Hunan Provincial Key Laboratory of Intelligent Information
Processing and Application for Hengyang Normal University (No.
IIPA18K06).

\normalsize

\end{document}